\documentclass{daj}

\usepackage
{
        amssymb,
        amsmath,
        amsthm,
        nicefrac,
        bookmark
}

\usepackage[all]{xy}

\dajAUTHORdetails{%
  title = {A Union of Euclidean Metric Spaces is Euclidean}, 
  author = {Konstantin Makarychev and Yury Makarychev},
  plaintextauthor = {Konstantin Makarychev, Yury Makarychev},
  keywords = {metric geometry, local-global properties, Lipschitz extension},
}   

\dajEDITORdetails{%
   year={2016},
   number={14},
   received={22 March 2016},   
   published={10 August 2016},  
   doi={10.19086/da.876},       
}   

\newtheorem{theorem}{Theorem}[section]
\newtheorem{lemma}[theorem]{Lemma}

\newtheorem{claim}[theorem]{Claim}
\newtheorem{corollary}[theorem]{Corollary}
\newtheorem{definition}[theorem]{Definition}

\newtheorem{conjecture}{Conjecture}
\newtheorem{question}[conjecture]{Question}

\DeclareMathOperator {\Ball}{Ball}
\newcommand {\set}   [1] {\left\{ #1 \right\}}
\newcommand {\bbR}    {\mathbb{R}}
\newcommand {\calL}    {{\cal{L}}}
\newcommand {\Exp}       {\mathbb{E}}
\newcommand {\EE}    [2] {\Exp_{#1}\left[#2\right]}

\begin{document}

\begin{frontmatter}[classification=text]

\title{A Union of Euclidean Metric Spaces is Euclidean} 

\author[konstantin]{Konstantin Makarychev}
\author[yury]{Yury Makarychev\thanks{Supported by  NSF awards CAREER CCF-1150062 and IIS-1302662.}}

\begin{abstract}
Suppose that a metric space $X$ is the union of two metric subspaces $A$ and $B$ that embed into Euclidean space with distortions $D_A$ and $D_B$, respectively.
We prove that then $X$ embeds into Euclidean space with a bounded distortion (namely, with distortion at most $7D_A D_B + 2(D_A+D_B)$).
Our result settles an open problem posed by  Naor.
Additionally, we present some corollaries and extensions of this result.
In particular, we introduce and study a new concept of an ``external bi-Lipschitz extension''.

In the end of the paper, we list a few related open problems.
\end{abstract}
\end{frontmatter}

\section{Introduction}
In this paper, we give an affirmative answer to the following question posed by Assaf Naor.
\begin{question}\label{q:Naor}
 Consider a metric space $X$ that is the union of two metric subspaces $A$ and $B$. Suppose that $A$ and $B$ embed into Euclidean space with distortions $D_A$ and $D_B$, respectively. Does $X$ embed into Euclidean space with bounded distortion?
\end{question}
\noindent We prove that the metric space $X$ embeds into $\ell_2$ with distortion at most $7 D_A D_B + 2(D_A + D_B)$.

\medskip
We note that the question is related to the recent research in theoretical computer science on local--global properties of metric spaces and
the power of lift--and--project relaxations for combinatorial optimization problems~\cite{ALNRRV,CMM,CMM-SA}.
One of the main goals of this research is to understand how constraints on relatively small subsets of a metric space affect its global properties.
In particular,
Arora, Lov\'asz, Newman, Rabani, Rabinovich, and Vempala~\cite{ALNRRV} asked the following question.
\begin{question} \label{q:ALNRRV}
For $n\geq k \geq 1$ and $p \in[1,\infty)$, find the least value of $D = D_{n,k,p}$ such that the following is true.
If $X$ is a finite metric space on $n$ points such that every $k$-point subset of $X$ isometrically embeds into $\ell_p$,
then $X$ embeds into $\ell_p$ with distortion at most $D_{n,k,p}$.
\end{question}
In~\cite{CMM}, Charikar, Makarychev, and Makarychev
showed that
$$c \left(\frac{\log n}{\log k + \log \log n}\right)^{1/p} \leq D_{n,k,p} \leq C \log (n/k),$$
for some positive absolute constants $c$ and $C$. They also showed that an upper bound of $C (\log (n/k) + \log\log (1/\alpha))$
holds if not all but only an $\alpha$ fraction of all subsets of size $k$ isometrically embed into $\ell_p$.
In contrast, in this paper, we assume only that \textit{two} subsets of $X$ (subsets $A$ and $B$) embed isometrically into $\ell_2$;
note that at least one of them should be of size $n/2$.

\begin{figure}[t]
\footnotesize
\centering
\begin{tabular}{|c|c|c|c|c|}
  \hline
    & \multicolumn{2}{c|}{ assumption on subsets of $X$ that embed into $\ell_p$} &  &    \\
   & number of subsets & size of each subset & distortion with which $X \hookrightarrow\ell_p$& value of $p$  \\ \hline
  Q\ref{q:Naor} & 2, subsets partition $X$ & avg. size is at least $n/2$ & $O(1)$ & $p=2$  \\ \hline
  Q\ref{q:ALNRRV} & an $ \alpha$ fraction of subsets & $k$ & $O(\log(n/k) + \log\log (1/\alpha)$ &  $p\in[1,\infty]$  \\
  \hline
\end{tabular}
\caption{Comparison of the assumptions and guarantees for Questions~\ref{q:Naor} and~\ref{q:ALNRRV}.}
\end{figure}

An analog of Question~\ref{q:Naor} for ultrametrics was studied by Mendel and Naor in~\cite{MN}. They showed
that if $(X, d)$ is the union of   $(A, d)$ and $(B,d)$
that embed into ultrametric spaces with distortions $D_A$ and $D_B$, then $X$ also embeds into an ultrametric space with distortion at most $(D_A+2)(D_B+2) -2$.

\subsection{Preliminaries}\label{sec:prelim}
We denote the $k$ dimensional Euclidean space by $\ell_2^k$ and the (separable) Hilbert space by $\ell_2^\infty$.
The Lipschitz constant of a map $f$ from a metric space $(X, d_X)$ to a metric space $(Y, d_Y)$ is
$\|f\|_{Lip} = \sup_{\substack{x,y\in X\\ x\neq y}} \frac{d_Y(f(x),f(y))}{d_X(x,y)}$. The distortion
of an embedding $f:X \hookrightarrow Y$ is $\|f\|_{Lip} \|f^{-1}\|_{Lip}$ (where $f^{-1}$ is the inverse map from
$f(X)\subset Y$ to $X$).

For a metric space $X$ and a Banach space $V$, let the Lipschitz extension constant $e_k(X, V)$ be the minimal constant $C$ such that the following holds:
for every subset $Y$ of $X$ of size at most $k$ and every map $f:Y \to V$ there exists an extension $\hat f: X \to V$ such that
$\|\hat f\|_{Lip} \leq C \|f\|_{Lip}$.
In this paper, we  use the Kirszbraun theorem that states that every map $f$ from a subset $Y$ of $\ell_2^a$ to $\ell_2^b$ can be extended
 to a map $\tilde f:\ell_2^a \to \ell_2^b$ so that $\|\tilde f\|_{Lip} = \|f\|_{Lip}$; in particular,  $e_k(\ell_2^a, \ell_2^b) = 1$ for every $k$~\cite{Kirszbraun} (see also~\cite{WW}).

\subsection{Our Results}
We prove the following theorem that answers Question~\ref{q:Naor}
affirmatively.
\begin{theorem}\label{thm:main}
Consider a metric spaces $(X,d)$. Assume that $X$ is the union of two metric subspaces  $A$ and $B$ that embed into $\ell_2^a$ and $\ell_2^b$ with distortions
$D_A$ and $D_B$, respectively. Then $X$ embeds into $\ell_2^{a+b+1}$ with distortion at most
$$7D_AD_B + 2(D_A+D_B).$$
If $D_A=D_B = 1$, then $X$ embeds into $\ell_2^{a+b+1}$ with distortion at most $8.93$.

\smallskip

 In this theorem, $a$ and $b$ may be finite or infinite.
\end{theorem}

As we show in~Lemma~\ref{lem:compactness} in Section~\ref{sec:compactness}, it is sufficient to prove Theorem~\ref{thm:main} only for finite metric spaces ---
the result for arbitrary metric spaces follows from the result for finite metric spaces.
So we assume below that $X$ is finite.

Let ${\cal D}(D_A, D_B)$ be the minimal $D$ such that every metric space $(X,d)$ as in Theorem~\ref{thm:main} embeds into $\ell_2$ with distortion at most $D$.
It is interesting to understand the dependence of ${\cal D}(D_A, D_B)$ on $D_A$ and $D_B$.
The result of Theorem~\ref{thm:main} can be restated as
$${\cal D}(D_A, D_B) \leq 7D_AD_B + 2(D_A+D_B) \quad\text{and}\quad {\cal D}(1,1) < 8.93\quad.$$
There is a trivial lower bound on ${\cal D}(D_A,D_B)$: ${\cal D}(D_A,D_B) \geq \max(D_A, D_B)$.  We also prove that
${\cal D}(D_A,D_B) \geq 3$.
\begin{theorem}\label{thm:lowerbound}
For every $\varepsilon > 0$, there exists a finite metric space $(X,d)$, which is
the union of two metric subspaces  $A$ and $B$ such that
\begin{itemize}
\item $A$ and $B$ embed into $\ell_2^{n-1}$ isometrically (where $n = |A| = |B|$),
\item any embedding of $X$ into $\ell_2\equiv \ell_2^{\infty}$ has distortion at least $3-\varepsilon$.
\end{itemize}
\end{theorem}

We present some corollaries and extensions to Theorem~\ref{thm:main} in Section~\ref{sec:extensions}. In particular, we
introduce and study a new notion of an ``external bi-Lipschitz extension''.

\section{Proof of Theorem~\ref{thm:main}}
In this section, we prove Theorem~\ref{thm:main}.
We now give a brief outline of our proof.
Since $A$ embeds into $\ell_2^a$ with distortion $D_A$, and $B$ embeds into $\ell_2^b$ with distortion $D_B$, there are non-contracting embeddings
$\varphi_A: A \hookrightarrow \ell_2^a$ and $\varphi_B: B \hookrightarrow \ell_2^b$ with $\|\varphi_A\|_{Lip} \leq D_A$ and $\|\varphi_B\|_{Lip} \leq D_B$.
Our first goal is to construct a map $\psi_B = \psi: X \to \ell_2^b$ that is Lipschitz on $X$ and bi-Lipschitz on $B$ (see Lemma~\ref{lem:psi}).
We start with proving Lemmas~\ref{lem:a-cover} and~\ref{lem:f-norm}. We define $\psi_B$
by letting $\psi_B = \varphi_B$ and then extending $\psi$ to a Lipschitz map from $X$ to $\ell_2^b$;
we use Lemmas~\ref{lem:a-cover} and~\ref{lem:f-norm} to show that the extension exists.
Similarly to $\psi_B$, we  construct a map $\psi_A$ that is Lipschitz on $X$ and bi-Lipschitz on $A$.
Then we consider the direct sum $\psi_A \oplus \psi_B$. This map is Lipschitz on $X$ and bi-Lipschitz on $A$ and on $B$;
however, it is not necessarily a bi-Lipschitz embedding of $X$ since it may significantly decrease distances between points  in $A$ and $B$.
Finally, we consider a map $\psi_{\Delta}:X \to {\mathbb R}$, which, loosely speaking, preserves distances between points in $A$ and $B$, and
obtain a desired embedding $\Psi = \psi_A \oplus \psi_B\oplus \psi_{\Delta}$.

\begin{definition}
For a point $x$ in a metric space $(X,d)$ and a radius $r\geq 0$, we denote the ball of radius $r$ around $x$ by
$\Ball_r(x) = \set{y:d(x,y) \leq r}$.
\end{definition}

\begin{definition}
For every $a\in A$, define $R_a = d(a, B)$; for every $b\in B$, define $R_b = d(b, A)$,
where $d(a,B)$ and $d(b,A)$ denote the distances from $a$ to the set $B$ and from $b$ to the set $A$,
respectively.
\end{definition}

\begin{definition}
Let $\alpha > 0$. We say that $A' \subset A$ is an $\alpha$-cover for $A$ with respect to $B$ if it satisfies the following two properties.
\begin{enumerate}
  \item For every $a\in A$, there is $a'\in A'$ such that $R_{a'} \leq R_a$ and $d(a,a') \leq \alpha R_a$.
  \item For every distinct $a_1', a_2' \in A'$, we have $d(a_1',a_2') \geq \alpha \min(R_{a_1'}, R_{a_2'})$.
\end{enumerate}
\end{definition}

\begin{lemma}\label{lem:a-cover}
Assume that $X= A \cup B$ is a finite metric space.
For every $\alpha > 0$, there exists an $\alpha$-cover $A'$ for $A$ with respect to $B$.
\end{lemma}
\begin{proof}
We prove the lemma by induction on the size of $A$. If $A=\varnothing$, we let $A' = \varnothing$; then $A'$ satisfies conditions (1) and (2).
Assume now that the statement of the lemma holds if $|A| < k$ (for some $k\geq 1$). We prove that it holds if $|A| = k$.
Find a point $u$ in $A$ with the smallest value of $R_{u}$; that is, a point $u \in A$ closest to $B$.
Let $Z = A \setminus \Ball_{\alpha R_{u}} (u)$. Note that $|Z| < |A| = k$. By the induction hypothesis, the statement of the lemma
holds for $Z$. Let $Z'$ be an $\alpha$-cover for $Z$ with respect to $B$, and $A' = Z' \cup \set{u}$.
We claim that $A'$ is an $\alpha$-cover for $A$ w.r.t.~$B$. We verify that $A'$ satisfies properties (1) and (2) of an $\alpha$-cover.

1. Let $a$ be a vertex in $A$. Consider two possibilities. Assume first that $a\in \Ball_{\alpha R_{u}} (u)$. Then let $a' = u$. By our choice of $u$,
$R_{a'} \leq R_a$. Since $a \in \Ball_{\alpha R_{u}} (u)$, we also have $d(a,a') = d(a,u)  \leq \alpha R_a$. Thus, property (1) holds.

Assume now that
$a\notin \Ball_{\alpha R_{u}} (u)$. Then $a\in Z$. By the induction hypothesis, there exists $a'\in Z' \subset A'$ such that $R_{a'} \leq R_a$ and $d(a,a') \leq \alpha R_a$, as required.

2. Consider $a_1', a_2'\in A'$. If both $a_1'$ and $a_2'$ are in $Z'$, then $d(a_1',a_2') \geq \alpha \min(R_{a_1'}, R_{a_2'})$ by the induction hypothesis.
So let us assume that either $a_1'$ or $a_2'$ is not in $Z'$. That is, $a_1' = u$ or $a_2' =u$. Without loss of generality, we assume that $a_1' = u$ and $a_2' \in Z'$.
Since $a_2' \in Z' \subset Z$, we have $a_2'\notin\Ball_{\alpha R_{u}}(u)$. Hence,
$d(a_1', a_2') \geq \alpha R_{u} = \alpha R_{a_1'} \geq \alpha \min(R_{a_1'}, R_{a_2'})$.
\end{proof}

Let $A'$ be an $\alpha$-cover.
Consider a map $f:A'\to B$ that maps every point $a'\in A'$ to a point in $B$ closest to $a'$ (we break ties arbitrarily).
That is, $f(a')$ is such that
$$d(a', f(a')) = d(a', B) = R_{a'}.$$
We show that $f$ is a Lipschitz map.
\begin{lemma}\label{lem:f-norm}
We have, $\|f\|_{Lip} \leq 2 (1 + 1 / \alpha)$.
\end{lemma}
\begin{proof}
Let $a_1'$ and $a_2'$ be two points in $A'$.
\begin{align*}
d(f(a_1'), f(a_2')) & \leq d(f(a_1'), a_1') + d(a_1', a_2') + d(a_2', f(a_2')) = R_{a_1'} + R_{a_2'} + d(a_1', a_2') \\
&= 2 \min( R_{a_1'}, R_{a_2'}) + | R_{a_1'} - R_{a_2'}| + d(a_1', a_2') .
\end{align*}
Note that $\min(R_{a_1'}, R_{a_2'}) \leq d(a_1',a_2')/\alpha$ by property 2 of an $\alpha$-cover, and
$| R_{a_1'} - R_{a_2'}| = |d(a_1', B) - d(a_2',B)| \leq d(a_1', a_2')$. Therefore,
$$d(f(a_1'), f(a_2')) \leq 2 d(a_1',a_2')/\alpha + d(a_1',a_2') + d(a_1',a_2') = 2(1+1/\alpha) \cdot d(a_1', a_2').$$
\end{proof}
By combining the maps $f$ and $\varphi_B$ we can obtain a
$2(1+1/\alpha) D_B$-Lipschitz embedding of $A'$ to $\ell_2^b$. We now show how to extend this embedding to
the entire set $X$.
\begin{lemma}\label{lem:psi}
Assume that $X$ is finite.
There exists a map $\psi:X \to \ell_2^b$ such that
\begin{enumerate}
\item For every $a_1,a_2\in A$,
\begin{equation}\label{eq:psi-a-a}
\|\psi(a_1) - \psi(a_2)\| \leq 2(1+1/\alpha) D_A D_B d(a_1,a_2).
\end{equation}
\item For every $b_1,b_2\in B$,
\begin{equation}\label{eq:psi-b-b}
d(b_1,b_2) \leq \|\psi(b_1) - \psi(b_2)\| = \|\varphi_B(b_1) - \varphi_B(b_2)\| \leq  D_B d(b_1,b_2).
\end{equation}
\item For every $a\in A$ and $b\in B$,
$$d(a,b) - (1+ \alpha) (2 D_A D_B +  1) R_a \leq \|\psi(a) - \psi(b)\| \leq (2(1+\alpha) D_A D_B + (2+\alpha) D_B) d(a,b).
$$
\end{enumerate}
\end{lemma}
\begin{proof}

\begin{figure}\label{fig:com-diag}
  \centering
$$
\xymatrix @C=6.7pc {
A \ar[d]^{\varphi_A}   & \ar[l]_{\subset}
A'
\ar[d]^{\varphi_A} \ar[r]^{f} &
B\ar[d]^{\varphi_B}\\
\varphi_A(A) \subset \ell_2^a \ar@/_1.8pc/[rr]^{\tilde g} &\ar[l]_{\subset}\varphi_A(A') \subset \ell_2^a \ar[r]^{ g = \varphi_B^{\vphantom{-1}} f \varphi_A^{-1}}  &
\ell_2^b
}
$$
\caption{Commutative diagram for maps $g$ and $\tilde{g}$.}
\end{figure}
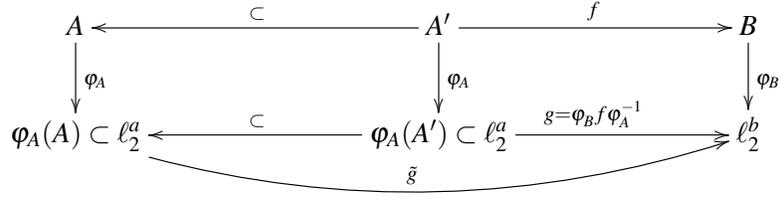
Consider map $g = \varphi_B^{\vphantom{-1}} f \varphi_A^{-1}$ from $\varphi_A(A')$ to $\varphi_B(B) \subset \ell_2^b$ (see Figure~\ref{fig:com-diag}).
We upper bound the Lipschitz   norm of $g$ (using Lemma~\ref{lem:f-norm})
$$\|g\|_{Lip} \leq  \|\varphi_B\|_{Lip} \|f\|_{Lip} \|\varphi_A^{-1}\|_{Lip} \leq D_B \cdot 2 (1 + 1/\alpha) \cdot 1.$$
By the Kirszbraun theorem, there is an extension $\tilde g:\ell_2^a\to \ell_2^b$ of $g$ to $\ell_2^a$ with $\|\tilde g\|_{Lip} = \|g\|_{Lip} \leq 2 (1 + \alpha) D_B/{\alpha}$.
Define map $\psi: X \to\ell_2^b$ as follows:
$$
\psi(x) =
\begin{cases}
\tilde g(\varphi_A(x)), &\text{if } x\in A,\\
\varphi_B(x), & \text{if } x\in B.
\end{cases}
$$
Note that if $x\in A\cap B$ then $x$ must be in $A'$ and thus $f(x) = x$. Therefore, both formulas for $\psi(x)$,  $\tilde g(\varphi_A(x))$ and $\varphi_B(x)$, are equal,
and $\psi(x)$ is well defined.
We prove that the map $\psi$ satisfies the conditions of the lemma.
\smallskip

1. Consider $a_1, a_2\in A$. We have,
$$\|\psi(a_1) - \psi(a_2)\| = \|\tilde g(\varphi_A(a_1)) - \tilde g(\varphi_A(a_2))\| \leq \|\tilde g\|_{Lip} \|\varphi_A\|_{Lip} \,d(a_1,a_2)
\leq2(1+1/\alpha) D_A D_B \cdot d(a_1,a_2).
$$

2. Consider $b_1,b_2\in B$. By the definition of $\psi$, $\|\psi(b_1) - \psi(b_2)\| = \|\varphi_B(b_1) - \varphi_B(b_2)\|$. Since $\varphi_B$ is a non-expanding map with Lipschitz constant $D_B$, we have
$$d(b_1,b_2) \leq \|\psi(b_1) - \psi(b_2)\| \leq  D_B d(b_1,b_2).$$

3. Finally, consider $a\in A$ and $b\in B$. By property 1 of an $\alpha$-cover, there is an $a'$ in $A'$ such that $R_{a'} \leq R_a$ and $d(a, a') \leq \alpha R_a$. Let $b' = f(a')$.
Note that $d(a',b') = R_{a'} \leq R_a$ and $\psi(a') = \varphi_B(f(a')) = \psi(b')$.
We have,
$$\|\psi(a) - \psi(b)\| \leq \|\psi(a) - \psi(a')\| + \|\psi(a') - \psi(b')\| + \|\psi(b') - \psi(b)\| = \|\psi(a) - \psi(a')\| + \|\psi(b') - \psi(b)\|.$$
Plugging in upper bounds (\ref{eq:psi-a-a}) and (\ref{eq:psi-b-b}) for $\|\psi(a) - \psi(a')\|$ and $\|\psi(b) - \psi(b')\|$, we get
\begin{equation}\label{eq:bound-psi}
\|\psi(a) - \psi(b)\| \leq 2(1+1/\alpha) D_A D_B d(a, a') + D_B d(b, b').
\end{equation}
We now upper bound $d(a,a')$ and $d(b,b')$. Note that $R_a = d(a, B) \leq d(a,b)$. Thus
\begin{equation}\label{eq:dist-a-aprime}
d(a,a') \leq \alpha R_a  \leq \alpha d(a, b),
\end{equation}
and
$$d(b,b') \leq d(b,a) + d(a,a') + d(a',b') \leq d(a,b) + \alpha d(a, b) + R_a \leq (2 + \alpha) d(a,b).$$
From (\ref{eq:bound-psi}), we get
$$
\|\psi(a) - \psi(b)\| \leq \bigl(2(1+\alpha) D_A D_B + (2 + \alpha) D_B\bigr) d(a,b)
$$

On the other hand, we have
$$\|\psi(a) - \psi(b)\| \geq \|\psi(b) - \psi(b')\| - \|\psi(a') - \psi(b')\| - \|\psi(a') - \psi(a)\| \geq d(b,b') - 0 -  2(1+ \alpha) D_A D_B R_a.$$
Here we used that $\psi(a') = \psi(b')$; we upper bounded $\|\psi(a') - \psi(a)\|$ using (\ref{eq:psi-a-a}) and the first inequality in~(\ref{eq:dist-a-aprime}).
We bound $d(b,b')$ as follows:
$$d(b,b') \geq d(a,b) - d(a,a') - d(a',b') \geq d(a,b) -\alpha R_a  - R_{a'} \geq d(a,b) - (\alpha + 1) R_a.$$
We get,
$$\|\psi(a) - \psi(b)\| \geq d(a,b) - (1+ \alpha) (2 D_A D_B +  1) R_a.$$
This concludes the proof of Lemma~\ref{lem:psi}.
\end{proof}

\noindent We are ready to prove Theorem~\ref{thm:main}.
\begin{theorem}[Theorem~\ref{thm:main} restated]\label{thm:main-restated}
Consider a metric spaces $(X,d)$. Assume that $X$ is the union of two metric subspaces  $A$ and $B$ that embed into $\ell_2^a$ and $\ell_2^b$ with distortions
$D_A$ and $D_B$, respectively. Then $X$ embeds into $\ell_2^{a+b+1}$ with distortion at most
$$7D_AD_B + 2(D_A+D_B).$$
If $D_A=D_B = 1$, then $X$ embeds into $\ell_2^{a+b+1}$ with distortion at most $8.93$.

 Furthermore, for given non-contracting embeddings $\varphi_A: A \hookrightarrow \ell_2^a$ and $\varphi_B: B \hookrightarrow \ell_2^b$ with $\|\varphi_A\|_{Lip} \leq D_A$ and $\|\varphi_B\|_{Lip} \leq D_B$, there is a non-contracting embedding $\Psi:X \hookrightarrow \ell_2^{a+b+1}$ with distortion at most
 $7D_AD_B + 2(D_A+D_B)$ such that $\|\Psi(u) - \Psi(v)\| \geq \|\varphi_A(u) - \varphi_A(v)\|$ for $u,v \in A$ and $\|\Psi(u) - \Psi(v)\| \geq \|\varphi_B(u) - \varphi_B(v)\|$ for $u,v \in B$.
 \smallskip

 In this theorem, $a$ and $b$ may be finite or infinite.
\end{theorem}
\begin{proof}
As we noted above, it is sufficient to prove the statement only for finite sets $X$;
see Lemma~\ref{lem:compactness} for details. So we assume that $X$ is finite.

We let $\alpha = 1/2$ in the proof of the bound $7D_AD_B + 2(D_A+D_B)$ on the distortion, which holds for every $D_A$ and $D_B$; we let $\alpha = 0.3114$ in the proof of a tighter bound of 8.93, which holds for $D_A=D_B=1$.

We construct a map $\psi:X\to\ell_2^b$ as in Lemma~\ref{lem:psi}. Denote $\psi_B = \psi$. Similarly (switching $A$ with $B$ in the statement of Lemma~\ref{lem:psi}), we construct a map $\psi_A:X\to\ell_2^a$.

Let $\beta = (1+ \alpha) (2 D_A D_B +  1)$ and $\gamma = \sqrt{1/2} \,\beta$. Define a map $\psi_{\Delta}: X \to \bbR$ by $\psi_{\Delta}(a) = \gamma R_a$ for $a\in A$ and
 $\psi_{\Delta}(b) = - \gamma R_b$ for $b\in B$. Finally, define an embedding $\Psi: X \to \ell_2^{a+b+1}$:
$$\Psi(x) = \psi_A(x) \oplus \psi_B(x) \oplus \psi_\Delta(x) \in \ell_2^{a+b+1}.$$
First, we show that $\Psi$ is a non-contracting map. For $a_1, a_2 \in A$,
$$\|\Psi(a_1) - \Psi(a_2)\| \geq \|\psi_A(a_1) - \psi_A(a_2)\| \stackrel{\text{\tiny by (\ref{eq:psi-b-b})}}{\geq} d(a_1, a_2).$$
Similarly, for $b_1, b_2 \in B$, $\|\Psi(b_1) - \Psi(b_2)\| \geq d(b_1, b_2)$.
For $a\in A$ and $b\in B$, we have
$$\|\Psi(a) - \Psi(b)\|^2 = \|\psi_A(a) - \psi_A(b)\|^2 + \|\psi_B(a) - \psi_B(b)\|^2 + \|\psi_{\Delta}(a) - \psi_{\Delta}(b)\|^2.$$
By Lemma~\ref{lem:psi}, item 3, and the definition of $\psi_{\Delta}$,
\begin{align*}
\|\psi_A(a) - \psi_A(b)\| &\geq d(a,b) - \beta R_b,\\
\|\psi_B(a) - \psi_B(b)\| &\geq d(a,b) - \beta R_a,\\
\|\psi_{\Delta}(a) - \psi_{\Delta}(b)\| &= \gamma (R_a + R_b). 
\end{align*}
By a simple case analysis, we show that these bounds imply that  $\|\psi(a) - \psi(b)\| \geq d(a,b)$.
\begin{claim}
$$\|\psi(a) - \psi(b)\| \geq d(a,b).$$
\end{claim}
\begin{proof}
We assume without loss of generality that $R_a \leq R_b$. Consider three cases. Assume first that $\beta R_b \leq d(a,b)$. Then
\begin{align*}
\|\psi(a) - \psi(b)\|^2 &\geq (d(a,b) - \beta R_a)^2 + (d(a,b) - \beta R_b)^2 + \beta^2 (R_a + R_b)^2/2\\
&=d(a,b)^2 + (d(a,b) - \beta R_a - \beta R_b)^2 + \beta^2 (R_a - R_b)^2 /2 \geq d(a,b)^2.
\end{align*}
Assume now that $\beta R_a \leq d(a,b) \leq \beta R_b$. Then
\begin{align*}
\|\psi(a) - \psi(b)\|^2 &\geq (d(a,b) - \beta R_a)^2 + \beta^2 (R_a + R_b)^2/2 \\
&= d(a,b)^2   -2 \beta d(a,b)  R_a + \beta^2 (3R_a^2 + 2 R_aR_b+ R_b^2)/2\\
&\geq d(a,b)^2   -2 \beta^2  R_aR_b + \beta^2 (3R_a^2 + 2 R_aR_b+ R_b^2)/2\\
&=d(a,b)^2 + \beta^2 ((\sqrt{3}R_a - R_b)^2 + 2(\sqrt{3} -1)R_aR_b)/2 \geq d(a,b)^{2}.
\end{align*}
Finally, assume that $d(a,b) \leq \beta R_a$. Then
$$\|\psi(a) - \psi(b)\|^2 \geq \beta^2 (R_a + R_b)^2/2 \geq  2 \beta^2 R_a^2 \geq 2d(a,b)^2.$$
\end{proof}
We conclude that $\Psi$ is a non-contracting map.
We now upper bound the Lipschitz constant of~$\Psi$.
For $a_1, a_2\in A$, we have
\begin{align*}
\|\Psi(a_1) - \Psi(a_2)\|^2 &= \|\psi_A(a_1) - \psi_A(a_2)\|^2 + \|\psi_B(a_1) - \psi_B(a_2)\|^2 + \|\psi_{\Delta}(a_1) - \psi_{\Delta}(a_2)\|^2 \\
&\leq (D_A^2 + 4 (1 + 1/\alpha)^2 D_A^2 D_B^2 ) d(a_1,a_2)^2 + \gamma^2 (R_{a_1}- R_{a_2})^2 \\
&\leq (D_A^2 + 4 (1 + 1/\alpha)^2 D_A^2 D_B^2 + \gamma^2) d(a_1,a_2)^2,
\end{align*}
here, we used that $|R_{a_1} - R_{a_2}| \leq d(a_1,a_2)$.
Similarly, for $b_1, b_2 \in B$, we have
$$\|\Psi(b_1) - \Psi(b_2)\|^2 \leq (D_B^2 + 4 (1 + 1/\alpha)^2 D_A^2 D_B^2 + \gamma^2) d(b_1,b_2)^2.$$
Now consider $a\in A$ and $b\in B$. Denote $\xi_A = 2 (1 + \alpha) D_A D_B + (2 + \alpha)D_A$ and $\xi_B = 2 (1 + \alpha) D_A D_B + (2 + \alpha)D_B$. We have,
\begin{align*}
\|\Psi(a) - \Psi(b)\|^2 &\leq  \|\psi_A(a) - \psi_A(b)\|^2 + \|\psi_B(a) - \psi_B(b)\|^2 + \|\psi_{\Delta}(a) - \psi_{\Delta}(b)\|^2
\\
&\leq
(\xi_A^2 + \xi_B^2) d(a,b)^2 + \gamma^2(R_a+R_b)^2 \leq (\xi_A^2 + \xi_B^2 + 4 \gamma^2\bigr) d(a,b)^2,
\end{align*}
here, we used that $R_a\leq d(a,b)$ and $R_b \leq d(a,b)$.

We first derive an upper bound for the Lipschitz constant of $\Psi$ in the general case (for arbitrary
$D_A$ and $D_B$). We plug in $\alpha = 1/2$ (this value of $\alpha$ is suboptimal, but we use it
to simplify the calculations), and get
\begin{align*}
\gamma^2 &= \frac{9}{8}(2D_A D_B + 1)^2,\\
\frac{\|\Psi(a_1)- \Psi(a_2)\|^2}{d(a_1,a_2)^{{2}}} &\leq \frac{81}{2} D_A^2 D_B^2 + \frac{9}{2} D_A D_B + D_A^2 +   \frac{9}{8}    \leq (7 D_A D_B + 2 D_A + 2  D_B)^2,\\
\frac{\|\Psi(b_1)- \Psi(b_2)\|^2}{d(b_1,b_2)^{{2}}} &\leq (7 D_A D_B + 2 D_A + 2  D_B)^2,\\
\frac{\|\Psi(a)- \Psi(b)\|^2}{d(a,b)^{{2}}} &\leq
36 D_A^2 D_B^2 + 15(D_A + D_B)D_A D_B + \frac{25}{4} (D_A^2 + D_B^2) + 18 D_A D_B + \frac{9}{2} \\
&\leq 36 D_A^2 D_B^2 +  15(D_A + D_B)D_A D_B + \frac{25}{4} (D_A^2 + D_B^2) + 18 D_A D_B + \frac{9}{2} \\&
\  {}
{}+\left(13 (D_A^2 D_B^2 + D_A^2 D_B + D_B^2 D_A) - \frac{9}{4} (D_A^2 + D_B^2) -10 D_A D_B - \frac{9}{2}\right)\\
&= (7 D_A D_B + 2 D_A + 2  D_B)^2.
\end{align*}
We get that the Lipschitz constant of $\Psi$ is at most $7D_AD_B + 2(D_A+D_B)$. Since $\Psi$ is non-contracting, the distortion of $\Psi$ is at most $7D_AD_B + 2(D_A+D_B)$.

Now consider the special case $D_A = D_B = 1$. Let $\alpha = 0.3114$.
\begin{align*}
\gamma^2 &= \frac{9(1+\alpha)^2}{2},\\
\frac{\|\Psi(a_1)- \Psi(a_2)\|^2}{d(a_1,a_2)^2} &\leq 1 + 4(1 + 1/\alpha)^2 + \frac{9(1+\alpha)^2}{2} < 79.7 < 8.93^2,\\
\frac{\|\Psi(b_1)- \Psi(b_2)\|^2}{d(b_1,b_2)^2} &\leq 1 + 4(1 + 1/\alpha)^2 + \frac{9(1+\alpha)^2}{2} < 79.7 < 8.93^2,\\
\frac{\|\Psi(a)- \Psi(b)\|^2}{d(a,b)^2} &\leq
2 (2 + \alpha + 2 (1 + \alpha))^2 + 18 (1 + \alpha)^2 <  79.7 < 8.93^2.
\end{align*}
We get that the distortion of $\Psi$ is less than $8.93$.
\end{proof}

\subsection{Compactness argument}~\label{sec:compactness}
In this section, we prove that it is sufficient to prove Theorem~\ref{thm:main-restated} only for finite metric spaces.
\begin{lemma}\label{lem:compactness}
I. Let $(X,d)$ be a separable metric space. Assume that every finite subset of $X$ embeds into $V = \ell_2^m$ with distortion at most $D$ (where $m$ is either finite or infinite).
Then $(X,d)$ embeds into $\ell_2^m$ with distortion at most $D$.

II. Furthermore, assume that $X = A \cup B$. Let $\alpha:A\times A\to {\mathbb R}$ and $\beta:B\times B\to {\mathbb R}$ be continuous functions.
Assume that for every finite subset $Y$ of $X$ there is an non-contracting embedding $f:Y\hookrightarrow V$ with $\|f\|_{Lip} \leq D$ such that
$\|f(x) - f(y) \| \leq \alpha(x,y)$ for every $x,y\in A\cap Y$ and $\|f(x) - f(y) \| \leq \beta(x,y)$ for every $x,y\in B \cap Y$.
Then there exists an non-contracting embedding $f:X \hookrightarrow V$ with $\|f\|_{Lip} \leq D$  such that
$\|f(x) - f(y) \| \leq \alpha(x,y)$ for every $x,y\in A$ and $\|f(x) - f(y) \| \leq \beta(x,y)$ for every $x,y\in B$.
\end{lemma}
\begin{proof}
I. Let $x_0, x_1,\dots, x_k, \dots$ be a dense sequence in $X$. Consider an orthonormal basis
$e_1,\dots, e_m$ in $\ell_2^m$, if $m < \infty$, or $e_1,\dots, e_k,\dots$ in $\ell_2^\infty$,
if $m = \infty$. Let $V_k$ be the linear span of vectors $e_1,\dots, e_k$ if $k \leq m$ and $V_k = V$ if $k > m$. In particular, let $V_0 = \set{0}$.

By the condition of the lemma, for every $k$, there exists an embedding $f_k$ of $\set{x_0,\dots, x_k}$ to $V$ with distortion at most $D$.
We  assume without loss of generality that
$$d(x_i, x_j) \leq \|f_k(x_i) - f_k(x_j)\| \leq D d(x_i, x_j).$$
Further, we assume that $f_k(x_0) = 0$. Applying the Gram--Schmidt process to
$f_k(x_0), \dots, f_k(x_k)$,
we get an isometry $T:V\to V$ such that $Tf_k(x_i) \in V_i$ for $i\in\set{0,\dots, k}$ (in particular, $f_k(x_0) = 0$).
Denote $f_k' = Tf_k$. Additionally, extend $f_k'$ to $\set{x_i: 0 \leq i < \infty}$ by letting $f_k'(x_i) = 0$ for $i > k$.

Let $B_i = \set{x\in V_i: \|x\| \leq D d(x_i, x_0)}$. Note that $f_k'(x_i) \in B_i$ since either
\begin{itemize}
\item $f_k'(x_i) = 0 \in B_i$, or
\item $\|f_k'(x_i)\| = \|f_k'(x_i) - f_k'(x_0)\| \leq D d(x_i,x_0)$ and $f_k'(x_i) \in V_i$, thus
$f_k'(x_i) \in B_i$.
\end{itemize}

Therefore, $f_k'\in \prod_{i=0}^\infty B_i$. Since each set $B_i$ is compact, the space $\prod_{i=0}^\infty B_i$ is compact;
also, the space is metrizable, since it is a countable product of metric spaces.
Therefore, the sequence $f_1',\dots, f_k', \dots$ has an accumulation point. Denote it by $f$.

Since for every $i$ and $j$, and all sufficiently large $k$, $d(x_i, x_j) \leq \|f_k'(x_i) - f_k'(x_j)\| \leq D d(x_i, x_j)$,
we have that $d(x_i, x_j) \leq \|f(x_i) - f(x_j)\| \leq D d(x_i, x_j)$. We obtained a desired embedding $f$.
Finally, we extend $f$ from $\set{x_i}$ to $X$ by continuity, and obtain an embedding of $X$ to $V$ with distortion at most $D$.

II. We proceed as in item I, except that we choose sequence $x_0, x_1,\dots, x_k, \dots$ so that $\set{x_i}\cap A$ is dense in $A$ and
$\set{x_i}\cap B$ is dense in $B$. Then since every function $f_k$ satisfies the additional requirements, we get by continuity
that $f$ also satisfies them.
\end{proof}

\section{Lower bound on distortion}\label{sec:lowerbound}
In this section, we prove Theorem~\ref{thm:lowerbound}.
Before we proceed to the proof, recall the definition and some properties of graph Laplacians.
Consider a graph $G = (V, E)$. The Laplacian $\calL$ of $G$ is a matrix with entries
$$\calL_{uv} =
\begin{cases}
\deg u,& \text{ if } u =v,\\
-1, &\text{ if } (u,v) \in E,\\
0, &\text{ otherwise.}
\end{cases}
$$
We write $X\succeq Y$ for two symmetric $n\times n$ matrices $X$ and $Y$, if $X-Y$ is positive semidefinite.

\begin{lemma} \label{lem:edge-decomp}
Let $G = (A \cup B, E)$ be a complete bipartite graph with parts $A$ and $B$ of size $n$ each.
There exists a partition of the set of edges $E$ into two disjoint sets $E_1$ and $E_2$ such that
$$\frac{1}{(1+\delta)} \calL_1 \preceq \frac{1}{2} \calL \preceq (1+\delta) \calL_1 \qquad\text{and}\qquad \frac{1}{(1+\delta)} \calL_2 \preceq \frac{1}{2} \calL \preceq (1+\delta) \calL_2,$$
where $\delta = c/\sqrt{n}$ (for some absolute constant $c$), and $\calL$, $\calL_1$ and $\calL_2$ are the Laplacians of $G$, $G_1= (A\cup B,E_1)$ and $G_2 = (A\cup B, E_2)$, respectively.
\end{lemma}
\begin{proof}
Let $E_1$ be a random subset of $E$ chosen uniformly among all subsets of $E$. Let $E_2=E\setminus E_1$. Then $G_1$ is a random bipartite $G(n,n,1/2)$ graph.
Thus $(1+\delta)^{-1} \calL_1 \leq \frac{1}{2} \calL \leq (1+\delta)\calL_1$ with probability at least $2/3$ (see \cite{FK}).
Similarly, $G_2$ is a random bipartite $G(n,n,1/2)$ graph, and
$(1+\delta)^{-1} \calL_2 \leq \frac{1}{2} \calL \leq (1+\delta)\calL_2$ with probability at least $2/3$
($G_1$ and $G_2$ are, of course, not independent).
Therefore, with probability at least $1/3$, $E_1$ and $E_2$ satisfy the conditions of the lemma.
\end{proof}
\begin{lemma}\label{lem:Pineq}
Let $G = (A \cup B, E)$, $E_1$ and $E_2$ be as in Lemma~\ref{lem:edge-decomp}.
Consider a map $f:A\cup B \to \ell_2$. Then for $i\in\set{1,2}$
$$\frac{1}{(1+\delta)^2} \, \EE{(u,v)\in E}{\|f(u) - f(v)\|^2} \leq \EE{(u,v)\in E_i}{\|f(u) - f(v)\|^2}
\leq (1 + \delta)^2 \, \EE{(u,v)\in E}{\|f(u) - f(v)\|^2}$$
\end{lemma}
\begin{proof}
Let $h:A\cup B\to \{0,1\}$ be a function equal to $0$ on $A$ and $1$ on $B$.
Then,
$$|E_i| = \langle h , \calL_i h\rangle \leq (1+\delta) \langle h , \calL h\rangle/2 =  (1+\delta) |E|/2. $$
Let $f_j(u)$ be the $j$-th coordinate of $f(u)$ in some fixed orthonormal basis. We have,
\begin{align*}
\Exp_{(u,v)\in E_i}& \bigl[ \|f(u) - f(v)\|^2\bigr] = \frac{1}{|E_i|} \sum_{(u,v)\in E_i}{ \|f(u) - f(v)\|^2} =
\frac{1}{|E_i|}\sum_j \langle f_j, \calL_i f_j\rangle \\
&\geq \frac{2}{(1+\delta)|E|} \cdot
 \sum_j \langle f_j, \Bigl(\frac{1}{2(1+\delta)} \calL\Bigr) f_j\rangle =
  \frac{1}{(1+\delta)^2} \cdot \frac{1}{|E|}\sum_{j} \langle f_j, \calL f_j\rangle
 \\
 & =  \frac{1}{(1+\delta)^2} \cdot \frac{\sum_{(u,v)\in E} \|f(u) - f(v)\|^2}{|E|}
  =
 \frac{1}{(1+\delta)^2} \cdot \EE{(u,v)\in E}{ \|f(u) - f(v)\|^2}.
 \end{align*}
 The proof of the other part of the inequality is analogous.
\end{proof}

Now we are ready to prove Theorem~\ref{thm:lowerbound}.
\begin{proof}[Proof of Theorem~~\ref{thm:lowerbound}]
Assume that $\varepsilon < 1$. Let $\delta = \varepsilon /6$. Choose $n$ so that $c/\sqrt{n} < \delta$.
Let $A$ and $B$ be two disjoint sets consisting of $n$ vertices each. Consider the complete bipartite graph $G=(A \cup B, E)$ with parts $A$ and $B$.
Partition all edges $E$ into two disjoint sets $E_1$ and $E_2$ as in Lemma~\ref{lem:edge-decomp}.
 Denote the Laplacians of $G= (A\cup B, E)$, $G_1=  (A\cup B, E_1)$ and $G_2= (A\cup B, E_2)$ by $\calL$, $\calL_1$ and $\calL_2$, respectively.

Define a metric space on $X = A \cup B$ as follows
$$d(u,v) =
\begin{cases}
0,&\text{if } u = v, \\
2,&\text{if } u\neq v \text{ and either both }a,b \in A \text{ or both } u,v \in B,\\
1,&\text{if } (u,v) \in E_1,\\
3,&\text{if } (u,v) \in E_2.
\end{cases}
$$
It is easy to see that $(X, d)$ is a metric space. Indeed, we have  $d(u,v) \geq 0$; $d(u,v) = 0$ if and only if $u = v$; $d(u,v) = d(v,u)$.
Consider three distinct vertices $u$, $v$, and $w$. If $u,w\in A$ or $u,w\in B$ then
$d(u,w) = 2 = 1 + 1\leq d(u,v) + d(v,w)$. If $u\in A$ and $w\in B$ then either $d(u,v) =2$ or $d(v,w) = 2$, and, therefore, $d(u,w)\leq 3 = 1 + 2 \leq d(u,v) + d(v,w)$.

Note that both metric spaces $(A, d)$ and $(B,d)$ are isometric to regular simplices in $\bbR^{n-1}$; in particular, they embed isometrically in $\bbR^{n-1}$.
We now show that every embedding of $(X,d)$ into $\ell_2$ has distortion at least $3 - \varepsilon$.
Consider an embedding $f$ of $(X,d)$ into $\ell_2$. Denote the distortion of $f$ by $D$.
 Then for some $c > 0$,
$$\|f(u) -f(v)\| \leq c D \text{ for } (u,v) \in E_1 \text{ and } \|f(u)-f(v)\| \geq 3c  \text{ for } (u,v) \in E_2.$$
Thus,
$$\frac{\EE{(u,v)\in E_1}{ \|f(u) - f(v)\|^2}}{\EE{(u,v)\in E_2}{ \|f(u) - f(v)\|^2}} \leq \frac{D^2}{9}.$$
On the other hand, by Lemma~\ref{lem:Pineq},
$$\frac{\EE{(u,v)\in E_1}{ \|f(u) - f(v)\|^2}}{\EE{(u,v)\in E_2}{ \|f(u) - f(v)\|^2}} \geq \left(\frac{1}{1+\delta}\right)^4.$$
Therefore, $D \geq 3/(1+\delta)^2 \geq 3 - \varepsilon$.
\end{proof}

\section{Corollaries and extensions to Theorem~\ref{thm:main}}\label{sec:extensions}
\subsection{Bi-Lipschitz extension}
In this section, we prove an ``external'' bi-Lipschitz extension theorem. Let $A\subset \ell_2^a$, $B\subset \ell_2^b$,
and $f$ be a Lipschitz map from $A$ to $B$. The Kirszbraun theorem states that the map $f$ can be extended to a Lipschitz
map $\tilde f$ from $\ell_2^a$ to $\ell_2^b$. Is there a counterpart of this theorem for bi-Lipschitz maps?
Note that there may be no bi-Lipschitz extension even for a map from a subset of $\mathbb R$ to a subset of $\mathbb R$.
Consider, for instance, a bi-Lipschitz map $f$ from $\set{0,1,2}\subset {\mathbb R}$ to $\set{0,1,2}\subset {\mathbb R}$ that maps $0$ to $0$, $1$ to $2$ and $2$ to $1$.
There is no continuous injective extension of $f$ to $\mathbb R$.

We prove, however, that there is an ``external'' extension of a bi-Lipschitz map.
\begin{definition}
Let $A\subset \ell_2^a$ and $B\subset \ell_2^b$. Let $f:A\to B$ be a bi-Lipschitz map. We say that a pair of maps $f_1:\ell_2^a \to \ell_2$ and $f_2:\ell_2^b\to \ell_2$ is an external extension
of $f$ with distortion $D$ if
\begin{itemize}
\item $f_1$ and $f_2$ have distortion at most $D$.
\item for every $a\in A$, $f_1(a) = f_2(f(a))$.
 \end{itemize}
\end{definition}
\begin{theorem}[External Bi-Lipschitz Extension Theorem]\label{thm:bi-Lip_ext}
 Let $A\subset U = \ell_2^a$ and $B\subset V= \ell_2^b$. Let $f:A\to B$ be a bi-Lipschitz map with distortion $D$.
 Then there exists an external extension $(f_1, f_2)$ of $f$ with distortion at most $O(D)$.
 \end{theorem}
 \begin{proof}
 Without loss of generality, we may assume that $f$ is non-contracting and $\|f\|_{Lip} \leq D$. We assume that $U$ and $V$ are disjoint
 (by replacing $V$ with an isometric copy of $V$ if necessary).
 We are now going to take the union of $U$ and $V$, identify each point $a\in A$ with $f(a) \in B$, and
 consider the shortest path metric on the obtained space.
 Formally, let $X = (U \cup V) / \set{(a,b): b = f(a)}$ (the union of $U$ and $V$ with each point $a\in A$ identified
 with $f(a)\in B$). Let $d$ be the shortest path metric on $X$; specifically, define metric $d$ as follows.
 \begin{itemize}
 \item For $u,v \in U$, let $d(u,v) = \|u - v\|$.
 \item For $u\in U$ and $v\in V$, let $d(u,v) = \inf_{x\in A} (\|u-x\| + \|f(x)- v\|)$.
 \item For $u,v \in V$, let
\begin{align*}d(u,v) &= \min (\|u-v\|, \inf_{w\in A} (d(u,w) + d(w,v))) \\
&= \min(\|u-v\|,\inf_{x,y\in B} (\|u-x\| + \|f^{-1}(x)-f^{-1}(y)\| + \|y-v\|)).
\end{align*}
 \end{itemize}
 Note that the distance between identified points $a\in A$ and $f(a)\in B$ is $0$. It is straightforward to check
  that $(X,d)$ is a metric space.

 We bound the distortions $D_A$ and $D_B$ with which $(U,d)$ and $(V,d)$ embed into $\ell_2^a$ and $\ell_2^b$, respectively.
 Note that $(U,d)$ is isometric to $(\ell_2^a,\|\cdot\|)$. So $D_A = 1$.
 Now consider the identity map $id$ from $(V, d)$ to $(\ell_2^b,\|\cdot\|)$.
 Note that $d(u,v) \leq \|u - v\|$ by the definition of $d$. So $id$ is non-contracting.
 Additionally, either $d(u,v) = \|u-v\|$ or
 \begin{align*}
 d(u,v) &= \inf_{x,y\in B} \|u-x\| + \|f^{-1}(x)-f^{-1}(y)\| + \|y-v\| \\
 &\geq \inf_{x,y\in B} \|u-x\| + \|x - y\|/D + \|y-v\| \geq \| u - v\| / D.
 \end{align*}
 So $\|id\|_{Lip} \leq D$.

 We apply\footnote{We would get better parameters if we applied Lemma~\ref{lem:psi} directly; however, for simplicity of
  exposition we use Theorem~\ref{thm:main-restated} here.} Theorem~\ref{thm:main-restated}
  to the metric space $X$ and get that $X$ embeds into $\ell_2$ with distortion at most $9 D + 2$.
 Denote the restrictions of this embedding to  $U$ and $V$ by $f_1$ and $f_2$, correspondingly.
We claim that $(f_1,f_2)$ is an external extension of $f$ with distortion at most $D' = 9 D + 2$.

 Indeed, for every $a\in A$, $f_2(f(a)) = f_1(a)$ since we identified $a$ and $f(a)$. Then $f_1$ has distortion at most $D'$.
 For $f_2$, we have
 \begin{align*}
 f_2(u,v) &\geq \|u-v\| &&\text{\footnotesize (by the ``furthermore'' clause in Theorem~\ref{thm:main-restated})}\\
 f_2(u,v) &\leq D' d(u,v) \leq D'\, \|u - v\| &&\mbox{\footnotesize (since  $f_2$ has distortion at most $D'$  as a map from $(X,d)$ to $\ell_2$)}
 \end{align*}
 We conclude that $f_2$ has distortion at most $D'$ as a map from $V$ to $\ell_2$.
\end{proof}

\subsection{An analog of Theorem~\ref{thm:main} for arbitrary normed spaces}
Theorem~\ref{thm:main} applies only to embeddings into Euclidean spaces:
if  $U$ and $V$ are \textit{Euclidean spaces}, $A$ and $B$ embed into $U$ and $V$ with distortions $D_A$ and $D_B$, then
$X = A \cup B$ embeds into $U \oplus V \oplus {\mathbb R}$ with distortion $O(D_A D_B)$.
However, observe that the proof of Theorem~\ref{thm:main}  uses only once that $U$ and $V$ are Euclidean --
when it extends map $g$ to a map $\tilde g$ using the Kirszbraun theorem.
In this section, we note that it is possible to  generalize Theorem~\ref{thm:main} to arbitrary normed spaces $U$ and $V$.
Since spaces $U$ and $V$ do not necessarily
satisfy the Kirszbraun theorem, our bound on the distortion with which $X$ embeds into $U \oplus V \oplus {\mathbb R}$ depends on the Lipschitz extension
constant for normed spaces $U$ and $V$.
\begin{theorem}\label{thm:extension}
Consider a metric spaces $(X,d)$. Assume that $X$ is the union of two finite metric subspaces  $A$ and $B$ that embed into normed spaces $U$ and $V$ with distortions
$D_A$ and $D_B$, respectively. Denote $E_A = e_{|A|}(A, V)$ and $E_B = e_{|B|}(B, U)$ (see Section~\ref{sec:prelim} for the definition of the Lipschitz
extension constants $e_k$).
Then $X$ embeds into $U \oplus V \oplus {\mathbb R}$ with distortion at most
$$O(E_A D_B + D_A E_B).$$
\end{theorem}
The proof of the theorem is almost identical to the proof of Theorem~\ref{thm:main}. The only difference is that
in the proof of Lemma~\ref{lem:psi}, we define $\psi(x)$ for $x\in A$ not as
$\tilde g(\varphi_A(x))$ but rather as an extension of map $\varphi_B f:A' \to V$ from $A'$ to $A$.
Then $\|\varphi_B f\| \leq D_B 2 (1 + 1/\alpha)$ and $\|\psi|_A\|_{Lip} \leq 2 E_A D_B (1 + 1/\alpha)$.

Now we present a corollary of  this theorem communicated to us by Naor.
\begin{corollary}[Naor]
Consider a metric space $X= A\cup B$ on $n$ points. Assume that $A$ embeds isometrically into normed space $U$ and $B$ into normed space $V$.
Then $X$ embeds into $U\oplus V \oplus {\mathbb R}$ with distortion at most $O(\log n/\log \log n)$.
\end{corollary}
\begin{proof}
We apply Theorem~\ref{thm:extension}. Note that $D_A = D_B = 1$. As Lee and Naor~\cite{LN} showed, the Lipschitz extension constant
$e_k(Y, W) = O(\log k /\log \log k)$ for every metric space $Y$ and normed space $W$. Therefore,
$E_A = O(\log n/\log \log n)$ and $E_B = O(\log n/\log \log n)$. The corollary follows.
\end{proof}

\section{Open problems}\label{sec:open}
In this section, we present several  open problems.
\begin{question}
Obtain tight or almost tight upper and lower bounds for ${\cal D} (D_A, D_B)$. In particular, find the asymptotic behavior
of ${\cal D}(D, D)$ when $D\to\infty$. Currently, we only know that $\Omega(D) \leq {\cal D}(D, D) \leq O(D^2)$. We conjecture that ${\cal D}(D, D) = \Theta(D^2)$.
\end{question}

\begin{question}
We know that ${\cal D}(1,1)\in [3,8.93)$. What is the exact value of ${\cal D}(1,1)$?
\end{question}

\begin{question}
Study the problem for spaces $\ell_p$ with $p\neq 2$. Assume that $X = A\cup B$, where $A$ and $B$ embed into $\ell_p$ isometrically.
Is it true that $X$ embeds into $\ell_p$ with bounded distortion? We conjecture that the answer is negative for every $p\notin\set{2,\infty}$.
\end{question}

\begin{question}
Assume that $X = A_1 \cup A_2 \cup \dots \cup A_k$, and each $A_i$ embeds isometrically into $\ell_2$. What is the smallest $D_k$ such that $X$
necessarily embeds into $\ell_2$ with distortion at most $D_k$. We can get from Theorem~\ref{thm:main} by induction on $k$ that $D_k \leq 2^{O(k)}$. On the other hand,
even if all sets $A_1,\dots, A_k$ are singletons (and thus trivially embed into $\ell_2$), $D_k \geq \Omega(\log k)$~\cite{LLR,Ostrovskii}.
\end{question}
Finally, we want to reiterate that Question~\ref{q:ALNRRV}  is still open (see the introduction).
Currently known lower and upper bounds for $D_{n,k,p}$ do not match. In particular, the following question is interesting.
\begin{question}
What is  the value of $D_{n,k,p}$ for $k = \sqrt{n}$ and $p=2$?
\end{question}

\section*{Acknowledgements}
We thank Assaf Naor for useful discussions.
We thank Tommaso Goldhirsch and an anonymous referee for their very valuable comments about the preliminary version of this paper,
and, in particular, for finding a typo in the statement of Lemma~\ref{lem:psi}, item 3, and pointing out that by fixing this typo, it is possible to get
a better bound on the distortion
in the isometric case when $D_A = D_B = 1$.


\begin{dajauthors}
\begin{authorinfo}[konstantin]
  Konstantin Makarychev\\
  Microsoft Research\\
  Redmond, WA, USA\\
  konstantin.makarychev\imageat{}gmail\imagedot{}com \\
  \url{http://konstantin.makarychev.net}
\end{authorinfo}
\begin{authorinfo}[yury]
  Yury Makarychev\\
  Toyota Technological Institute at Chicago (TTIC)\\
  Chicago, IL, USA\\
  yury\imageat{}ttic\imagedot{}edu \\
  \url{http://ttic.uchicago.edu/~yury}
\end{authorinfo}
\end{dajauthors}

\begin{thebibliography}{99}
\bibitem{ALNRRV}
S.~Arora, L.~Lov\'asz, I.~Newman, Y.~Rabani, Y.~Rabinovich, and S.~Vempala.
\newblock Local versus global properties of metric spaces.
\newblock In Proc.~of the Symposium on Discrete Algorithms, pp.~41--50, 2006.

\bibitem{CMM}
M.~Charikar, K.~Makarychev, and Y.~Makarychev.
\newblock Local global tradeoffs in metric embeddings.
\newblock SIAM Journal on Computing 39, no. 6 (2010): 2487--2512.

\bibitem{CMM-SA}
M.~Charikar, K.~Makarychev, and Y.~Makarychev.
\newblock Integrality gaps for Sherali--Adams relaxations.
\newblock In Proc.~of the Symposium on Theory of Computing, pp.~283--292, 2009.

\bibitem{FK} Z.~F{\"u}redi and J.~Koml{\'o}s.
\newblock The eigenvalues of random symmetric matrices.
\newblock Combinatorica 1, no. 3 (1981): 233--241.

\bibitem{Kirszbraun} M.~Kirszbraun.
\newblock {\"U}ber die zusammenziehende und Lipschitzsche Transformationen.
\newblock Fundamenta Mathematicae 22, no 1 (1934): 77--108.

\bibitem{LLR} N.~Linial, E.~London, and Y.~Rabinovich.
\newblock The geometry of graphs and some of its algorithmic applications.
\newblock Combinatorica 15, no. 2 (1995): 215--245.

\bibitem{LN} J.~R.~Lee and A.~Naor.
\newblock Absolute Lipschitz extendability.
\newblock Comptes Rendus Mathematique 338, no. 11 (2004): 859--862.

\bibitem{MN} M.~Mendel and A.~Naor.
\newblock Ultrametric skeletons.
\newblock Proc. Natl. Acad. Sci. 110, no. 48 (2013): 19256--19262.

\bibitem{Ostrovskii} M.~I.~Ostrovskii.
\newblock Metric embeddings: bilipschitz and coarse embeddings into Banach spaces.
\newblock Vol. 49. Walter de Gruyter, 2013.

\bibitem{WW} J.~H.~Wells and L.~R.~Williams.
\newblock Embeddings and extensions in analysis. Vol. 84. Springer Science \& Business Media, 2012.

\end{thebibliography}
\end{document}